\documentclass[10pt]{article}
\usepackage{amssymb}

\usepackage{mathrsfs}
\usepackage{threeparttable}
\usepackage{amsfonts}
\usepackage{amsmath}
\usepackage{amssymb}
\usepackage{enumerate}
\usepackage{hyperref}
\usepackage{latexsym}
\usepackage{xcolor}

\usepackage{amsmath,amssymb}\usepackage{cite}\usepackage{mathrsfs}\usepackage[amsmath,thmmarks]{ntheorem}
\usepackage{listings}
\usepackage[titletoc]{appendix}\allowdisplaybreaks
\usepackage{enumitem}
\usepackage{enumerate}
\newlist{thmenum}{enumerate}{1} 
\setlist[thmenum]{label=(\arabic*)}
\makeatletter

\newcommand{\Rmnum}[1]{\expandafter\@slowromancap\romannumeral  #1@}
\makeatother

\makeatletter
\renewcommand\theequation{\thesection.\@arabic\c@equation}
\makeatother

\newtheorem{thm}{ Theorem}[section]%
\newtheorem{lem}[thm]{ Lemma}%
\input cyracc.def

\def\D{\Delta}  

\def\PSL{\hbox{\rm PSL}}

\def\Cay{\hbox{\rm Cay }}

\def\qed{\hfill $\Box$}

\newcommand\ZZ{\mathrm C}
\def\Z{{\rm Z}}

\def\Q{{\rm Q}}
\def\D{{\rm D}}

\def\C{{\rm C}}

\def\N{{\rm N}}

\newcommand{\pf}{\noindent {\bf Proof. \ }}
\begin{document}
\title{ A Note on Subgroup Perfect Codes in Cayley Graphs }
%
\author {Jingjian Li$^a$, Binbin Li$^a$ and Xianglin Liu$^b$\thanks{Corresponding author at: College of General Education, Guangxi Vocational University of Agriculture. People's Republic of China.}\\
a, School of Mathematics and Information Science \& Guangxi Base, \\
   Tianyuan Mathematical Center in Southwest China  \\
\& Guangxi Center for Mathematical Research   \\
\& Center for Applied Mathematics of Guangxi(Guangxi University), \\
Guangxi University, \\
Nanning, Guangxi 530004, P. R. China.\\
b, College of General Education, Guangxi Vocational University of Agriculture,\\
Nanning, Guangxi 530007, P. R. China.
}
\date{}
\maketitle
\renewcommand{\thefootnote}{}
\footnotetext{
E-mail addresses: lijjhx@gxu.edu.cn(J.J. Li), libb@st.gxu.edu.cn(B.B. Li),yx331229641@163.com(X.X. Liu).}

\baselineskip=16pt

\vskip0.5cm

{\bf Abstract}. In this paper, we give a necessary and sufficient condition for a subgroup to be a perfect code for finite groups. As an application, we determine all subgroup perfect codes of  extraspecial $2$-groups and finite groups whose  Sylow $2$-subgroup is extraspecial.

{\bf Keywords}: Cayley graph, subgroup perfect codes, extraspecial $2$-group.

Mathematics Subject Classification: 05C25, 05C69, 94B25\\

\section {Introduction}
All groups considered are finite, and all graphs considered are finite, undirected and simple.
Let  $G$ be a group and $S$ be an inverse-closed subset of $G \setminus\{1\}$. The Cayley graph $\Cay(G, S)$ on $G$ with connection set $S$ is the graph with vertex set $G$ such that $x, y \in G$ are adjacent if and only if $yx^{-1}\in S$.
A subset $C$ of the vertex set of a graph $\Gamma $ is called a perfect code \cite{j}  in $\Gamma$ if $C$ is an independent set and every vertex of $\Gamma$ is at distance no more than $1$ to exactly one vertex of $C$, in particular,
a perfect code is also called an independent perfect dominating set \cite{l} or  efficient dominating set \cite{dej}. A subgroup $H$ of a group $G$ is called a perfect code of $G$ if there exists a Cayley graph $\Cay(G,S)$ of $G$ admiting $H$ as a perfect code.

In the research of perfect codes in graphs, especially perfect codes in Cayley graphs have received special attention\cite{c,ch,dej,den,f,h,kh,ku,l,m,o,z,z2,z3}.
For a finite group $G$, the question whether a subgroup of $G$ is a perfect code of $G$ has attracted a lot of attention.
In \cite{h}, Huang, Xia  and Zhou introduced the concept of subgroup perfect codes in Cayley graphs and gave a necessary and sufficient condition for a normal subgroup of $G$ to be a perfect code and determined all subgroup perfect codes of cyclic groups and dihedral groups.
In \cite{m}, Ma, Walls, Wang and Zhou  proposed that a group is called code-perfect if every subgroup of $G$ is a perfect code, and showed that a group is code-perfect if and only if it has no element of order $4$, also they determined all subgroup perfect codes of the generalized quaternion group.
More recently, Zhang\cite{z} gave a necessary and sufficient condition for  subgroup  perfect codes  by  characterising its Sylow $2$-subgroup and obtained that all maximal subgroups of the special linear group $\PSL(2,q)$ are its perfect codes.
Furthermore, many interesting results about the subgroup perfect codes in Cayley graph have been given, for example\cite{ch,kh,z2,z3}.

Let $G$ be an arbitrary finite group and $N\trianglelefteq G$. A non-trivial element $x\in G$ is called a square if there exists $y\in G$ such that $x=y^{2}$. Set $\Omega_1({G})= \{ g\in G : g^{2}=1 \}$ and $G_2$ is a Sylow $2$-subgroup of $G$. In particular, we denote $\Omega_1(G)N=\bigcup_{g\in \Omega_1(G)} Ng$.

In this paper, by characterising $\Omega_1(\N_G(H_2))$ and $\Omega_1(\N_G(H_2)/H_2)$,  we obtain a necessary and sufficient condition for subgroup $H$ to be a perfect code of $G$.
As an application, we determine all subgroup perfect codes of extraspecial $2$-groups and  the finite  groups  whose Sylow $2$-subgroups are extraspecial.
Our main results as follows:

\begin{thm} \label{theorem1.1}
Let $G$ be a group and $H\leqslant G$. Assume that $P$ is a Sylow $2$-subgroup of $G$ such that $\N_G(H_2)_{2}\leqslant P$. Then the following statements are equivalent:
\begin{itemize}[label=\arabic*,font=\upshape]
\item [(1)] $H_2$ is a perfect code of $P$;

\item [(2)] $\Omega_1(\N_G(H_2)_2/H_2)=\Omega_1(\N_G(H_2)_2)H_2$;

\item [(3)] $\Omega_1(\N_G(H_2)/H_2)=\Omega_1(\N_G(H_2))H_2$;

\item [(4)] $H$ is a perfect code of $G$.
\end{itemize}
\end{thm}

A finite 2-group $G$ is called extraspecial $2$-group if it's center $\Z(G)\cong\Z_2$ and $G/\Z(G)$ is elementary abelian. It is clear that $\exp G=4$ in this case.  Denote $G_{m, 1}$  the center product of $m$ dihedral groups with order $8$ for $m$ positive integer.

\begin{thm}\label{theorem1.2}
Let $G$ be an extraspecial $2$-group and $H\leqslant G$. Then $H$ is a perfect code of $G$ if and only if either $H$ is non-abelian, or $H$ is abelian with $H\ntrianglelefteq G$, or $H$ is a maximal abelian subgroup of $G\cong G_{m, 1}$.
\end{thm}

\begin{thm} \label{theorem1.3}
Let $G$ be a group with Sylow $2$-subgroup being extraspecial and $H\leqslant G$. Then $H$ is a perfect code of $G$ if and only if either $|H|$ is odd or
 one of the following holds:
\begin{itemize}[label=\arabic*,font=\upshape]
\item [(1)]  $H_2$ is non-abelian;

\item [(2)]  $H_2$ is abelian and $|\N_G(H_2)_{2}|<|G_2|$;

\item [(3)]  $H_{2}$ is abelian, $|H_2|^{2}=2|G_2|$ and $G_2\cong G_{m,1}$.
\end{itemize}
\end{thm}

At the end of this section we give some of the notation used in this paper. For group-theoretic terminology not defined here, we refer the reader to \cite{hu}.
For a group $G$, use $G'$ and $\Phi(G)$ to denote the commutator subgroup of $G$ and $\Phi(G)$ to denote the Frattini subgroup of $G$, respectively.
For an element $x$, denote the order of $x$ by $o(x)$.
The central product of $H$ and $K$ is denoted by $H \circ K$ for two subgroups $H$ and $K$.

\section{Preliminaries}
In this chapter, we will introduce some results to be used later in the form of technical lemma.

Let $H\leqslant G$ and $T\subseteq G$. Then $T$ is called a right transversal of $H$ in $G$ if $T$ contains exactly one element of each right coset of $H$, in particular, $T$ is called inverse-closed if $T=T^{-1}$.


\begin{lem}\cite[Theorem 1.2]{ch} \label{lemma2.1}
Let $G$ be a group and $H\leqslant G$. Then the following cases are equivalent:
\begin{itemize}[label=\arabic*,font=\upshape]
\item [(1)]  $H$ is a perfect code of $G$;

\item [(2)]  there exists an inverse-closed right transversal of $H$ in $G$;

\item [(3)]  for each $x \in G$ such that $x^{2}\in  H$ and $|H|/|H \cap H^{x}|$ is odd, there exists $y \in Hx$ such that $y^{2} = 1$;

\item [(4)]  for each $x\in G$ such that $HxH = Hx^{-1}H$ and $|H|/|H \cap H^{x}|$ is odd, there exists $y \in Hx$ such that $y^{2} = 1$.
\end{itemize}
\end{lem}

Note that $\Omega_1(G)= \{ g\in G : g^2=1 \}$. By considering the 2-subgroup or normal subgroup of a group, we can deduce that the following useful result.
\begin{lem} \label{lemma2.2}
Let $G$ be a group and $H\leqslant G$. Suppose that $H$ is  a $2$-group or $H\unlhd G$. Then $H$ is a perfect code of $G$ if and only if $\Omega_1(\N_G(H)/H)=\Omega_1(\N_G(H))H$.
\end{lem}
\pf Assume that $H$ is a perfect code of $G$. Then $H$ is also a perfect code of $\N_G(H)$.
Recall from  the case (2) of Lemma \ref{lemma2.1} that, $H$ has an inverse-closed right transversal $T$ of $H$ in $\N_G(H)$ such that $\N_G(H)=\bigcup_{t\in T}Ht$.
For any  $Hg\in \Omega_1(\N_G(H)/H)$, the definition of $\Omega_1(\N_G(H)/H)$ shows that  $Hg=(Hg)^{-1}$.
Since  $\N_G(H)/H=\bigcup_{\tiny{t\in T}} Ht$,  it follows that there exists $t\in T$ such that $Hg=Ht$ and hence $Ht=(Ht)^{-1}$.
Applying the definition of $T$ and $T=T^{-1}$, we see that $t$ is an involution, and so $Hg=Ht \in \Omega_1(\N_G(H))H$.
%
Moreover, since $\Omega_1(\N_G(H))H\subseteq \Omega_1(\N_G(H)/H)$, this implies  $\Omega_1(\N_G(H)/H)=\Omega_1(\N_G(H))H$.

Now suppose that $\Omega_1(\N_G(H)/H)=\Omega_1(\N_G(H))H$.
Let $x\in G$ such that $x^2\in H$ and $|H:H\cap H^x|$ is odd.
If $H$ is normal subgroup of $G$, then $\Omega_1(G/H)=\Omega_1(G)H$ and thus there exists an involution $y\in \Omega_1(G)$ such that $Hy=Hx$.
Deriving from the case (2) of Lemma \ref{lemma2.1} that $H$ is a perfect code of $G$.
Now assume that $H$ is $2$-subgroup. Suppose for a contradiction that $H$ is not a perfect code of $G$.
Applying the case (3) of Lemma \ref{lemma2.1},  there exists  $x \in G$ satisfying that $x^2\in H$,  $|H|/|H \cap H^x|$ is odd and $Hx$ has no involution. It follows from $H$ is a 2-subgroup that $x\in \N_G(H)$ and hence $Hx\in \Omega_1(\N_G(H)/H)$.
Since $\Omega_1(\N_G(H)/H)=\Omega_1(\N_G(H))H$,  there exists an involution $y\in\Omega_1(\N_G(H))$ such that $Hx=Hy$, which contradicts that $Hx$ has no involution. Therefore, $H$ is a perfect code of $G$.
$\Box$

Note that every elementary abelian $2$-group can be viewed as a vector space over $GF(2)$, and for an extraspecial $2$-group $G$, $\Z(G)\cong\ZZ_2$ and $G/\Z(G)$ is elementary abelian.
\begin{lem}\cite[\Rmnum{3}, Satz 13.7]{hu} \label{lemma2.3}
Let $G$ be an extraspecial $2$-group and $\Z(G)=\langle c \rangle\cong\ZZ_2$.
Define a bilinear function say $\beta$, from $G/\Z(G) \times G/\Z(G)$ to $GF(2)$  satisfies $\beta(\overline{x},\overline{y})=\alpha$, where $[x,y]=c^\alpha$, $\overline{x}=x\Z(G)$ and  $\overline{y}=y\Z(G)$.
Then $(G/\Z(G),\beta)$ is a non-degenerate symplectic space over  $GF(2)$, and $(\overline{x},\overline{y})$ is a hyperbolic pair of $(G/\Z(G),\beta)$ if and only if $xy\neq yx$.
\end{lem}

Following the Lemma \ref{lemma2.3}, we find that $G/\Z(G)$ is a non-degenerate symplectic space on $GF(2)$ and so the dimension of $G/\Z(G)$ is even. In this paper, we always denote by $2m$ the dimension of $G/\Z(G)$, where $m$ is a positive integer, and so the number of vectors in $G/Z(G)$ is $2^{2m}$. Leaded from the assumption that $\Z(G)\cong \ZZ_2$, we yield that the order of $G$ is $2^{2m+1}$.

 Note that a non-trivial element $g$ of group $G$ is called a square if there exists $h\in G$ such that $g=h^{2}$. The following lemma is about extraspecial $2$-group.

\begin{lem}\label{lemma2.4}
Let $G$ be an extraspecial $2$-group with order $2^{2m+1}$ and $H\leqslant G$. Then
\begin{itemize}[label=\arabic*,font=\upshape]
 \item [(1)]  $\Phi(G)=Z(G)=G'$ is of order $2$ and $G/Z(G)$ is elementary abelian $2$-group;
 \item [(2)] $G$ is isomorphic to one of the following:
  \begin{itemize}[label=\arabic*,font=\upshape]
    \item [({\romannumeral1})] the central product $G_{m,1}$ of $m$ dihedral groups;
    \item [({\romannumeral2})] the central product  $G_{m,2}$ of $m-1$ dihedral groups and quaternion group;
  \end{itemize}
 \item [(3)] every maximal abelian subgroup of $G_{m,1}$ is isomorphic to $\ZZ_4\times \ZZ_2^{m-1}$ or $\ZZ_2^{m+1}$ and every maximal abelian subgroup of $G_{m,2}$
is isomorphic to $\ZZ_{4}\times\ZZ_2^{m-1}$;
 \item [(4)]  $G$  has a unique  square element denoted by $c$, and $Z(G)=\langle c  \rangle$;
 \item [(5)]  $H \unlhd G$ if $\exp (H)=4$.
\end{itemize}
\end{lem}

\pf For the cases (1),  (2) and (3), we can be straightly get from the \cite[\Rmnum{3}, Satz 13.8]{hu}. In the following, we will prove the correctness of the cases (4) and (5).

Followed from case (1), we can deduce that $\exp G=4$.
Let $c$ be a non-trivial square element of $G$.
The \cite[\Rmnum{3}, Satz 3.14]{hu} shows that all of square elements of $G$ are contained in $\Phi(G)$, and so $c$ is contained in $\Phi(G)$. Again by the case (1), we yield that $\ZZ_2\cong\Phi(G)=\langle c \rangle$, that is $G$ has a unique non-trivial square element $c$ and hence the case (4) holds.

Assume that  $\exp (H)=4$.
Then $H$ has an element $x$ with order $4$ and $x^2$ is a non-trivial square element of $H$. Employing the case (4), one deduces that $x^2=c$ and $\Phi(G)=\langle c\rangle=\langle x^2\rangle\leqslant H$. Derived from $H/\Phi(G)\unlhd G/\Phi(G)$ as $G/\Phi(G)$ is elementary abelian, one yields that $H\unlhd G$ and hence the case (5) holds. $\Box$

Let $G$ be a finite group, $A\leqslant G$ and $B\leqslant G$. Then $G$ is called a central product of $A$ and $B$ denoted by $G=A \circ B$ if the following conditions hold: $G=AB$, $A\cap B\leqslant \Z(G)$ and $ab=ba$ for any two elements $a\in A$ and $b\in B$. The following Lemma \ref{lemma2.5} is simple but is useful which can be obtained from the proof of Proposition (a) of \cite[\Rmnum{3}, Satz 13.8]{hu}.
\begin{lem} \label{lemma2.5}
$\Q_8\circ\Q_8=\D_8\circ\D_8$.
\end{lem}

\section{Proof of Theorem \ref{theorem1.1}}
In this section, we will prove the correctness of Theorems \ref{theorem1.1}.\

\pf Let $G$ be a group and $H\leqslant G$. Assume that $P$ is a Sylow $2$-subgroup of $G$ such that $\N_G(H_2)_{2}\leqslant P$.

(1)$\Leftrightarrow$(2): Since $\N_{P}(H_2)=P\cap\N_G(H_2)=\N_G(H_2)_2$,  applying the Lemma \ref{lemma2.2} by replacing $G$ with $P$, we deduce that (1) and (2) are equivalent.

(2)$\Rightarrow$ (3): The definition of $\Omega_1(*)$ shows that $\Omega_1(\N_G(H_2))H_2\subseteq \Omega_1(\N_G(H_2)/H_2)$, where $*$ denotes an arbitrary finite group. In the following, we will prove that $\Omega_1(\N_G(H_2)/H_2)\subseteq \Omega_1(\N_G(H_2))H_2$.  Assume that (2) holds, i.e., $\Omega_1(\N_G(H_2)_2/H_2)=\Omega_1(\N_G(H_2)_2)H_2$.
Let $H_2x$ be an arbitrary element of $\Omega_1(\N_G(H_2)/H_2)$.
Then $(H_2x)^2=H_2$ and $x^2\in H_2$. It follows that $x$ is an even order element of $\N_G(H_2)$, and so there must exist an odd integer denoted by $n$, such that $x^n$ is a $2$-element of $\N_G(H_2)$ and $(x^n)^a\in\N_G(H_2)_2$ for $a\in\N_G(H_2)$.
Set $n=2k+1$ and $g=x^n$ for $k$ positive integer.
Since $H_2x$ is an involution of $\N_G(H_2)/H_2$, one yields that $H_2g=H_2x^n=(H_2x)^n=(H_2x)^{2k}H_2x=H_2x$, that is $H_2g$ is an involution.
Noting that $g^a=(x^n)^a\in \N_G(H_2)_2$ for $a\in\N_G(H_2)$ and $\Omega_1(\N_G(H_2)_2/H_2)=\Omega_1(\N_G(H_2)_2)H_2$, one deduces that there exists $y\in \Omega_1(\N_G(H_2))_2$ such that $H_2g^a=H_2y$, and so $H_2g=H_2y^{a^{-1}}$ as $H_2g^a=(H_2g)^a$. It follows that $H_2x=H_2g=H_2y^{a^{-1}}\in\Omega_1(\N_G(H_2))H_2$ and hence $\Omega_1(\N_G(H_2))H_2=\Omega_1(\N_G(H_2)/H_2)$ meeting expectations.

(3) $\Rightarrow$ (4):
Assume that the case (3) holds, that is $\Omega_1(\N_G(H_2)/H_2)=\Omega_1(\N_G(H_2))H_2$.
Suppose that $H$ is not a perfect code of $G$.
According to  Lemma \ref{lemma2.1} (3),  there exists  $g\in G\setminus H$ with $g^{2}\in H$  such that  $|H:H\cap H^{g}|$ is odd and  $Hg$ has no involution. Derived from $g\not\in H$ and $g^2\in H$, one yields that the order of $g$ is even.
Set $C=H\cap H^{g}$.
Since $g^{2}\in H$, one yields that $g^2=(g^2)^g\in H^g$ and so $g^{2}\in C$, in particular, $C^{g}=C$.
Let $\Delta$ be the set of all Sylow $2$-subgroups of $C$.
Then the Sylow-theorem shows that $|\Delta|$ is odd and considering the action of $\langle g\rangle$ by conjugation on $\Delta$, one yields that $\langle g\rangle$ must stabilize one element of $\Delta$ denoted by $D$ as the order of $g$ is even, that is $D^g=D$ and $g\in \N_G(D)$.
Applying the assumption that $|H:C|$ is odd, one yields that $D$ is also a Sylow 2-subgroup of $H$. It follows from the case (3) that $\Omega_1(\N_G(D)/D)=(\Omega_1(\N_G(D))D)$ and so $Dg=Dy$ for some
$y\in \Omega_1(\N_G(D))$. It follows that $Dg\subseteq Hg$ has an involution $y$, a contradiction.
Therefore, $H$ is a perfect code of $G$.


(4) $\Rightarrow$ (2):
Assume that (4) holds, that is $H$ is a subgroup perfect code of $G$. For any $H_2x\in \Omega_1(\N_G(H_2)_2/H_2)$ and $x\not\in H_2$,
one yields that $(H_2x)^2=H_2$, i.e., $H_2x^2=H_2$ and $x^2\in H_2$. Since $x\in\N_G(H_2)_2$ and $x^2\in H_2$, $(H\cap H^x)^x=H\cap H^x$ and $H_2\in H\cap H^x$, that is, $|H:H\cap H^{x}|$ is odd.
Since $H$ is a perfect code of $G$, by Lemma \ref{lemma2.1} (3),  there exists an involution $hx$ in $Hx$, for some $h\in H$.
Set $D=\langle H_2, x, h \rangle$. We claim that $\langle H_2,  h \rangle\unlhd D$ and so $D=\langle H_2,  h \rangle\langle x \rangle$. Since $x^{-1} h^{-1} x=x^{-1}h^{-1}((hx)^{-1}hx)x=(hx)^{-2}hx^{2}=hx^{2}\in\langle H_2, h \rangle$, we yield that $\langle x\rangle$ normalizes $\langle H_2, h \rangle$ as $x\in\N_G(H_2)$,  and hence the claim holds.  Derived from $x^2\in H_2\subseteq H$, one yields that $D=\langle H_2,  h \rangle\langle x \rangle=\langle H_2,  h \rangle\cup \langle H_2,  h \rangle x$, i.e., $|D:\langle H_2,  h \rangle|=2$. Noting that $x\in\N_G(H_2)_2$, $\langle H_2, h\rangle\leqslant H$, $x\not\in H_2$ and $H_2$ is a Sylow $2$-subgroup of $H$,
we obtain that $ H_2 \langle x \rangle$ is a Sylow $2$-subgroup of $D$.
Since $hx$ is an involution and $hx\in D$, let $y=hx$, one yields that $y$ is an involution of $D$. Then there exist some $h'x^{i}\in D$ such that $y^{h'x^{i}}\in H_2\langle x \rangle$, where $h'\in \langle H_2,  h \rangle$ and $i$ is a non-negative integer. It follows that $y^{h'}\in (H_2\langle x \rangle)^{x^{-i}}=H_2\langle x \rangle$.
Since $x\notin H$ and $h\in H$,  we deduce that $y^{h'}=(hx)^{h'}\notin H$ and so  $y^{h'}\notin H_2$.
Since $|H_2\langle x\rangle|=|H|o(x)/o(x^2)$ as $x^2\in H$, one yields that $|H_2\langle x\rangle:  H_2|=2$ and $H_2x=H_2y^{h'}\in \Omega_1(\N_G(H_2)_2)H_2$ and then  $\Omega_1 (\N_G(H_2)_2/H_2)=\Omega_1(\N_G(H_2)_2)H_2$.

To sum up, we have proved Theorem \ref{theorem1.1}.
\qed

\section{ Proof of Theorem \ref{theorem1.2} }
In this section, we will prove the Theorem \ref{theorem1.2} by some technical lemmas.
In this Section, we always let $G$ be an extraspecial 2-group with order $2^{2m+1}$ and $H\leqslant G$, where $m$ is a positive integer. Then $\Z(G)\cong\ZZ_2$ and $G/\Z(G)$ is elementary abelian, in particular, $\exp G=4$.

\begin{lem}\label{non-abelian} If $H$ is a non-abelian subgroup of $G$, then $H$ is a perfect code of $G$.
\end{lem}
\pf Since $H$ is not abelian, $\exp(H)\neq 2$ and there exist $x, y\in H$ such that $xy\neq yx$. Let $A=\langle x,y \rangle$. In the light of the fact that $\exp G=4$, we conclude that $\exp A=4$ and $o(x)$ and $o(y)$ are both the elements of set $\{2, 4\}$. The case (5) of Lemma \ref{lemma2.4} shows that $\langle x, y\rangle=A\unlhd G$.
Assume that $o(x)=o(y)=2$. Then $A$ is a dihedral group, and following from $\exp A=4$, one yields that $A\cong\D_8$.
Suppose that one of the $o(x)$ and $o(y)$ is 2, and the other is 4. Without loss of generality, we let $o(x)=4$ and $o(y)=2$. Applying the case (5) of Lemma \ref{lemma2.4}, one yields that $A\unlhd G$ and $\langle x\rangle\unlhd G$, in particular, $\langle x\rangle\unlhd A$. Since $xy\neq yx$, we deduce that $\langle x\rangle\cap \langle y\rangle=1$ and $A=\langle x\rangle\rtimes\langle y\rangle\cong\ZZ_4\rtimes\ZZ_2\cong\D_8$. Now that $o(x)=o(y)=4$. Then $\exp \langle x\rangle=\exp \langle y\rangle=\exp A=4$. Again by the case (5) of Lemma \ref{lemma2.4}, one deduces that $\langle x\rangle\unlhd G$ and $\langle y\rangle\unlhd G$. The case (4) of Lemma \ref{lemma2.4} shows that $G$  has a unique  square element and so $\langle x\rangle\cap \langle y\rangle=\Z(G)\cong\ZZ_2$, that is $\langle x,y \rangle\cong\Q_8$.

The Lemma \ref{lemma2.4} shows that $G'=\Z(G)=\langle c \rangle\cong \ZZ_2$ and so $[x, y]=c$ as $xy\neq yx$. Applying the Lemma \ref{lemma2.3}, one yields that there exists a bilinear function say $\beta$, from $G/\Z(G)\times G/\Z(G)$ to $GF(2)$ such that $(G/\Z(G),\beta)$ is a non-degenerate symplectic space over $GF(2)$ with
$(x\Z(G),y\Z(G))$ being a hyperbolic pair of $(G/\Z(G), \beta)$ as $xy\neq yx$(see the definition of $\beta$ given in the Lemma \ref{lemma2.3}).
Note that $G/\Z(G) = \langle x\Z(G),y\Z(G)\rangle \perp \langle x\Z(G),y\Z(G)\rangle^{\perp}$.
 Now we consider the bijection from the subgroups of $G/\Z(G)$ to the subspace of the symplectic space $(G/\Z(G),\beta)$. Set $B/\Z(G)$ is a preimage of $\langle x\Z(G),y\Z(G)\rangle^{\perp}$ in $G/\Z(G)$.
Since $Z(G)=\langle c\rangle=[x, y]\leqslant \langle x,y \rangle$, we deduce that the preimage of $\langle x\Z(G),y\Z(G)\rangle$ is $\langle x, y\rangle/\Z(G)$ and $G/\Z(G)=\langle x,y \rangle/\Z(G)B/\Z(G)$.
On the other hand, for any $d\in B$, $\beta(x\Z(G),d\Z(G))=\beta(y\Z(G),d\Z(G))=0$, combing with the definition of $\beta$, we can deduce that $[x,d]=[y,d]=1$ and $d\in \C_G(\langle x, y\rangle)$, that is $B\leqslant C_G(\langle x,y \rangle)$.
Therefore, $G/\Z(G)=\langle x,y \rangle/\Z(G)B/\Z(G)=\langle x,y \rangle B/\Z(G)=\langle x,y \rangle \C_G(\langle x,y \rangle)/\Z(G)=A\C_G(A)/\Z(G)$ and $G=A\C_G(A)$.

Noticing that $\exp A=4$, $\exp G=4$ and $A\leqslant H\leqslant G$, one yields that $\exp H=4$ and $H\trianglelefteq G$ (see the case (5) of  Lemma \ref{lemma2.4} for example).
Let $Hg\in \Omega_1(G/H)$ such that $o(g)=4$ and $g=hz$ for $h\in A$ and $z\in \C_G(A)$ as $G=A\C_G(A)$.
The case (4) of the Lemma \ref{lemma2.4} shows that there is a unique square element in $G$, say $c$. We claim that one of $o(h)$ and $o(z)$ is $2$, and another is $4$. If $o(h)=o(z)=4$, then $hz=zh$ leads to that $(hz)^2=h^2z^2=cc=1$, a contradiction arisen and hence the claim holds.
Assume that $o(h)=4$ and $o(z)=2$. Then $Hz=Hhz=Hg\in\Omega_1(G)H$ as $h\in A\leqslant H$.
Suppose that $o(h)=2$ and $o(z)=4$.
Note that $A\cong\D_8$ or $Q_8$. Take an element $h_1\in A$ such that $o(h_1)=4$. Then $h_1z=zh_1$ and $(h_1z)^2=(h_1)^{2}z^{2}=cc=1$, i.e., $h_1z$ is an involution of $G$. It follows that $Hg=Hhz=Hz=H(h_1z)\in\Omega_1(G)H$ and so $\Omega_1(G/H)\subseteq\Omega_1(G)H$. It is clear that $\Omega_1(G)H\subseteq\Omega_1(G/H)$ and hence $\Omega_1(G/H)=\Omega_1(G)H$. Further, applying the result of the Theorem \ref{theorem1.1}, we deduce that $H$ is a perfect code of $G$ and so the Lemma \ref{non-abelian} holds.\qed

\begin{lem}\label{nonormal} If $H\ntrianglelefteq G$, then $H$ is a  perfect code of $G$.
\end{lem}
\pf The assumption shows that $G$ is an extraspecial 2-group, $\exp G=4$ and $H\leqslant G$. Then $\exp H=2$ or $4$. Since $H\ntrianglelefteq G$, applying the case (5) of the Lemma \ref{lemma2.4}, one deduces that $\exp H=2$, i.e., $H$ is an elementary abelian $2$-subgroup of $G$.
Suppose that $H$ has a square element of $G$. Applying the case (4) of the Lemma \ref{lemma2.4}, $G$ has a unique square element denoted by $c$ and so $\Z(G)=\langle c\rangle\in H$ and $G/\Z(G)$ is abelian, which leads to that $H/\Z(G)\unlhd G/\Z(G)$ and $H\unlhd G$, contradicting to the assumption that $H\ntrianglelefteq G$.
Thus, $H$ has no square element of $G$. On the other hand, take an arbitrary element of $\Omega_1(\N_G(H)/H)$ denoted by $Hg$, then $g^2\in H$ leads to that $g^2=1$, that is $Hg\in\Omega_1(\N_G(H))H$ and $\Omega_1(\N_G(H)/H)\subseteq\Omega_1(\N_G(H))H$. The result of natural establishment $\Omega_1(\N_G(H))H\subseteq\Omega_1(\N_G(H)/H)$ shows that $\Omega_1(\N_G(H)/H)=\Omega_1(\N_G(H))H$. Applying the Theorem \ref{theorem1.1} we deduce that $H$ is a perfect code of $G$ and hence the Lemma \ref{nonormal} holds.\qed

\begin{lem} \label{normal non-maximal}If $H$ is a normal non-maximal abelian subgroup in $G$, then $H$ is not a perfect code of $G$.
\end{lem}
\pf Since $H$ is normal in $G$, by \cite[\Rmnum{3}  Satz 2.6]{hu}, we deduce that $H\cap\Z(G)\neq 1$ and further, $\Z(G)\leqslant H$ as $\Z(G)\cong\ZZ_2$.
 Noting that $H$ is a non-maximal abelian subgroup of $G$,  there must exist a maximal abelian subgroup of $G$ denoted by $A$ such that $H<A$. It follows that $H\unlhd A$ and $\Z(G)\unlhd A$.

Let $\overline{G}=G/\Z(G)$, $\overline{H}=H/\Z(G)$ and $\overline{A}=A/\Z(G)$. According to Lemma \ref{lemma2.3}, there exists a bilinear function $\beta$ such that $(\overline{G},\beta)$ is a non-degenerate symplectic space over $GF(2)$. Without causing confusion, we also simple denote $\overline{G}$ a space. It is clear that $\overline{H}$ and $\overline{A}$ are both subspace of $\overline{G}$.
Let $\overline{x_1}$, $\overline{x_2}$, $\cdots$, $\overline{x_r}$ be a basis of vector space $\overline{H}$ for $r$ positive integer.
Then $\overline{H}=\langle \overline{x_1}, \overline{x_2}, \cdots, \overline{x_r} \rangle={ \overline{\langle x_1, x_2, \cdots, x_r  \rangle\Z(G)} }$.
Therefore, $H=\langle x_1, x_2, \cdots, x_r  \rangle\Z(G)$.
Since $A$ is a maximal abelian subgroup of $G$ and $|G|=2^{2m+1}$, applying the case (3) of the Lemma \ref{lemma2.4}, we conclude that $|A|=2^{m+1}$.
It follows that the number of vectors of subspace $\overline{A}$ is $2^{m}$ and so the dimension of subspace $\overline{A}$ is $m$.
Now we extend the above basis of $\overline{H}$ to a basis of $\overline{A}$, say $\overline{x_1} $, $\overline{x_2} $, $\cdots$, $\overline{x_r}$, $\overline{x_{r+1}}$, $\cdots$,  $\overline{x_m}$.
Noting that $[x_i,x_j]=1$ for $1\leqslant i,j\leqslant m$, combining with the definition of $\beta$, we deduce that  $\beta(\overline{x_i}, \overline{x_j})=0$.
 Since the dimension of $\overline{A}$ and $\overline{G}$ are $m$ and $2m$ respectively, we conclude that $\overline{A}=\langle \overline{x_1}, \cdots, \overline{x_m}\rangle$ is a maximal totally isotropic subspace of $\overline{G}$.

Since $\overline{G}$ is a non-degenerate symplectic space and $\beta(\overline{x_l}, \overline{x_j})=0$ for $1\leqslant l, j\leqslant m$, there exist vectors  $\overline{y_1} $, $\overline{y_2} $, $\cdots$, $\overline{y_m}$ such that $\overline{x_1}, \cdots, \overline{x_m}$,  $\overline{y_1}, \cdots, \overline{y_m}$ is a symplectic basis of $\overline{G}$ and $( \overline{x_i}, \overline{y_i})$ are hyperbolic pairs of $\overline{G}$, where $i=1,2, \cdots, m$.
Then  $\overline{G}=\langle \overline{x_1}, \overline{y_1} \rangle \perp \langle \overline{x_2}, \overline{y_2} \rangle \perp \cdots \perp \langle \overline{x_m}, \overline{y_m} \rangle$ and $\beta(\overline{x_i}, \overline{y_i})=1$.
By the definition of $\beta$, we obtain that $[x_i,y_i]=c$, where $\langle c\rangle=Z(G)$.

Noting that $\langle \overline{x_i}, \overline{y_i} \rangle \subseteq \langle \overline{x_j}, \overline{y_j}\rangle^{\perp}$ for $i\neq j$, one yields that $\beta(\overline{x_i},\overline{x_j})=\beta(\overline{x_i},\overline{y_j})=\beta(\overline{y_i},\overline{x_j})=\beta(\overline{y_i},\overline{y_j})=0.$ By the definition of $\beta$, it concludes that $[x_i,x_j]=[x_i,y_j]=[y_i,x_j]=[y_i,y_j]=1$ and so
$\langle x_i, y_i \rangle\leqslant \C_G(\langle x_j, y_j \rangle)$.
Since $\Z(G)\leqslant \langle x_i,y_i \rangle$ and  $\overline{G}=\langle \overline{x_1}, \overline{y_1},  \overline{x_2}, \overline{y_2}, \cdots,  \overline{x_m}, \overline{y_m} \rangle$,
we deduce that $\overline{\langle x_i,y_i \rangle}=\langle x_i,y_i \rangle/\Z(G)$ and
$\overline{G}=\overline{\langle x_1, y_1 \rangle  \langle x_2, y_2 \rangle \cdots \langle x_m, y_m\rangle}$.
Therefore,  $G=\langle x_1, y_1 \rangle  \langle x_2, y_2 \rangle \cdots \langle x_m, y_m\rangle$.
In light of the fact that $\Z(G)\leqslant \langle x_i, y_i\rangle$ and $\Z(G)\leqslant \langle x_j, y_j\rangle$ such that $\Z(G)\leqslant \langle x_i, y_i\rangle \cap \langle x_j, y_j\rangle$.
If $\Z(G)< \langle x_i, y_i\rangle \cap \langle x_j, y_j\rangle$, then there exists $g\in (\langle x_i, y_i\rangle \cap \langle x_j, y_j\rangle)\setminus \Z(G)$ such that $\overline{g}\in \langle \overline{x_i}, \overline{y_i}\rangle \cap \langle \overline{x_j}, \overline{y_j}\rangle$ contradicting with the fact that $\langle \overline{x_i}, \overline{y_i}\rangle \cap \langle \overline{x_j}, \overline{y_j}\rangle=\overline{1}$. Thus $\Z(G)= \langle x_i, y_i\rangle \cap \langle x_j, y_j\rangle$, for $i\neq j$.
It follows that $ \langle x_i, y_i\rangle  \langle x_j, y_j\rangle= \langle x_i, y_i\rangle \circ \langle x_j, y_j\rangle$ for $i\neq j$.
For any $k$  such that $k\neq i$ and $k\neq j$, since $\langle \overline{x_i}, \overline{y_i}, \overline{x_j}, \overline{y_j}\rangle\subseteq \langle \overline{x_k}, \overline{y_k} \rangle^{\perp}$, by the same argument as above, we obtain that $\Z(G)=\langle x_i, y_i, x_j, y_j \rangle \cap \langle x_k, y_k \rangle$ and  $\langle x_i, y_i, x_j, y_j \rangle \langle x_k, y_k \rangle=\langle x_i, y_i, x_j, y_j \rangle \circ \langle x_k, y_k \rangle$.
By analogy, we conclude that
$G=\langle x_1, y_1 \rangle \circ \langle x_2, y_2 \rangle \circ \cdots \circ \langle x_m, y_m\rangle$.

Since  $H$ is not a maximal abelian subgroup of $G$, one yields that $r<m$, that is $x_{r+1}y_{r+1}\neq y_{r+1}x_{r+1}$. Now substituting $x_{r+1}$ and $y_{r+1}$ for $x$ and $y$ in the first paragraph of the proof of the Lemma \ref{non-abelian} respectively, we can deduce that there exists an element $g_{r+1}\in \langle x_{r+1}, y_{r+1} \rangle$ such that $o(g_{r+1})=4$.
Let $I=\{ i : o(x_i)=4, 1\leq i\leq r  \}$ and $M=\{1,2,\cdots,r\}$.
For any $i\in I$, if $o(y_i)=2$, then we set $y_i'=y_i$, otherwise, set $y_i'=g_{r+1}y_i$. The uniqueness of the square element of $G$ leads to the fact that $y_i'$ is an involution in any case.
Again by the first paragraph of the proof of the Lemma \ref{non-abelian}, we conclude that $\langle x_i, y_i' \rangle\cong \D_8$ and $x_i y_i'$ is an involution for $i\in I$.
Let $g=\prod\limits_{i\in I} y_i' g_{r+1}$. Applying the case (4) of Lemma \ref{lemma2.4}, the unique square element $c$ of $G$ is in $\Z(G)\leqslant H$ and so $g^2\in H$.
 Note that $H=\langle x_1, x_2, \cdots, x_r \rangle Z(G)$. Let  $h=x_1^{a_1}x_2^{a_2}\cdots x_r^{a_r}c^{a_{r+1}}$ be an arbitrary element of $H$ for $1\leqslant i\leqslant r+1$ and $a_i$ non-negative integer.
 In light of the fact that $g_{r+1}\in\C_G(H)$ and $y_i'x_j=x_jy_i'$ for $i\neq j$, we deduce that $hg=\prod\limits_{l\in M\setminus I} x_l^{a_l} \cdot  \prod\limits_{l\in I} x_l^{a_l}y_l' \cdot g_{r+1}\cdot c^{a_{r+1}}$.
Noticing that the order of $\prod\limits_{j\in M\setminus I} x_j^{a_j}$, $\prod\limits_{i\in I} x_i^{a_i}y_i'$ and $c^{a_{r+1}}$ are at most $2$ and $o(g_{r+1})=4$, we conclude that $o(hg)=4$.
It follows from the arbitrariness of $h$  that $Hg$ has no involution. Derived from  $\Omega_1(G)H$ is a subset of $\Omega_1(G/H)$ and $g^2\in H$, one yields  that $Hg\in\Omega_1(G/H)\setminus\Omega_1(G)H$ and so $\Omega_1(G)H$ is a proper subset of $\Omega_1(G/H)$. Employing the result of the Theorem \ref{theorem1.1}, we can deduce that $H$ is not a perfect code of $G$ and hence the Lemma \ref{normal non-maximal} holds.           \qed

\begin{lem}\label{maximal abelian subgroup}  Let $H$ be a maximal abelian subgroup of $G$. Then $H$ is a perfect code of $G$ if and only if  $G\cong G_{m,1}$, the center product of $m$ dihedral groups $\D_8$.
\end{lem}
\pf Suppose that $H$ is a  maximal abelian subgroup of $G$ and $G$ is an extraspecial 2-group. Then $\Z(G)\leqslant H$ and $H/\Z(G)\leqslant G/\Z(G)$ which is abelian, and so $H/\Z(G)\unlhd G/\Z(G)$ and $H\unlhd G$. Applying the case (3) of Lemma \ref{lemma2.4}, we find that $H$ is isomorphic to $\ZZ_4\times \ZZ_2^{m-1}$ or $\ZZ_2^{m+1}$, in particular, $|H|=2^{m+1}$.
According to the second paragraph in the proof of Lemma \ref{normal non-maximal},
 we can set $H=\langle x_1,x_2,\cdots, x_m\rangle \Z(G)$.
Furthermore, by the complete similar argument to the third and fourth paragraphs in the proof of Lemma \ref{normal non-maximal},
it is shown  that there exist $y_1, y_2,\cdots, y_m $ such that $G= \langle x_1, y_1  \rangle \circ \langle x_2, y_2  \rangle \circ \cdots  \circ  \langle x_m, y_m  \rangle$.
Let $M=\{1, 2,  \cdots, m\}$. From the first paragraph in the proof of Lemma \ref{non-abelian}, we conclude that  $\langle x_i, y_i \rangle\cong\D_8$ or $\Q_8$ for $i\in M$. Since $\langle x_i, y_i \rangle\langle x_j, y_j \rangle=\langle x_j, y_j \rangle\langle x_i, y_i \rangle$ for all $i$ and $j$ in $M$,
we can rearrange the factors such that for  some integer $t$,
$\langle x_i, y_i \rangle\cong \D_8$ for $1\leqslant i\leqslant t$  and  $\langle x_l, y_l \rangle\cong \Q_8$ for $t<l\leqslant m$.
Then $x_l$ and $y_l$ are of order $4$, where $t+1\leqslant l\leqslant m$.


\textbf{Claim 1:} If $G\cong G_{m,1}$, then $H$ is a perfect code of $G$.


Suppose that $G\cong G_{m,1}$. Applying the Lemma \ref{lemma2.5}, we conclude that $m-t$ is even.
In this paragraph, we always let $t+1\leqslant j\leqslant m-1$.
Let $z_j=y_jx_{j+1}$ and $w_{j+1}=y_{j+1}x_j$.
Further, by the case (4) of Lemma \ref{lemma2.4},  we deduce that $x_j^{2}=y_{j}^{2}=x_{j+1}^{2}=y_{j+1}^{2}=c$ for $\langle c \rangle=\Z(G)\cong\ZZ_2$ and $\Z(G)\leqslant\langle x_i, y_i\rangle$.
On the other hand, derived from $\langle x_j, y_j \rangle\leqslant \C_G(\langle x_{j+1}, y_{j+1} \rangle)$, one yields that
$z_j$ and $w_{j+1}$ are involutions.
 Further, combining with $\langle x_j, y_j \rangle\cong \langle x_{j+1}, y_{j+1} \rangle \cong \Q_8$, we conclude that $x_j^{z_j}=x_j^{y_jx_{j+1}}=x_j^{y_j}=x_j^{-1}$ and $x_{j+1}^{w_{j+1}}=x_{j+1}^{y_{j+1}x_j}=x_{j+1}^{y_{j+1}}=x_{j+1}^{-1}$ which leads to that $\langle x_j,z_j  \rangle\cong \langle x_{j+1},w_{j+1}  \rangle \cong\D_8$.
Simple calculation shows that
$$ z_{j}w_{j+1}
 =y_jx_{j+1} y_{j+1}x_j
 = y_j x_j  x_{j+1}y_{j+1}
 = y_j x_j x_j^{2}x_{j+1}y_{j+1}y_{j+1}^{2}=y_j x_j^3x_{j+1}y_{j+1}^3
 = y_j x_j^{-1}x_{j+1}y_{j+1}^{-1}.$$
On the other hand, $x_j^{y_j}=x_j^{-1}$ and $x_{j+1}^{y_{j+1}}=x_{j+1}^{-1}$ leads to $x_jy_j=y_jx_j^{-1}$ and $y_{j+1}x_{j+1}=x_{j+1}y_{j+1}^{-1}$ respectively, which follows that
$$ w_{j+1}z_{j}
  =y_{j+1}x_jy_jx_{j+1}
  =x_jy_jy_{j+1}x_{j+1}
  =y_j x_j^{-1}x_{j+1}y_{j+1}^{-1}=z_{j}w_{j+1}. $$
Therefore,  $\langle x_j, z_j \rangle\leqslant \C_G(\langle x_{j+1}, w_{j+1} \rangle)$.
By the complete similar argument to the fourth paragraph in the proof of Lemma \ref{normal non-maximal}, we deduce that $\Z(G)=\langle x_j, z_j \rangle\cap \langle x_{j+1}, w_{j+1}\rangle$. Then
$\langle x_j, y_j \rangle \circ\langle x_{j+1}, y_{j+1}\rangle=\langle x_j, y_j, x_{j+1}, y_{j+1} \rangle=\langle x_j,z_j, x_{j+1},w_{j+1}\rangle=\langle x_j, z_j \rangle \circ \langle x_{j+1},w_{j+1} \rangle$.
%
Derived from $G=\langle x_1, y_1  \rangle \circ \langle x_2, y_2  \rangle \circ \cdots  \circ  \langle x_m, y_m  \rangle$ and $m-t$ is even, we conclude that %
$$G=\langle x_1, y_1  \rangle \circ \langle x_2, y_2  \rangle \circ \cdots  \circ  \langle x_t, y_t  \rangle \circ \langle x_{t+1}, z_{t+1}  \rangle \circ \langle x_{t+2}, w_{j+2}  \rangle \circ
\cdots  \circ \langle x_{m-1}, z_{m-1}   \rangle  \circ \langle x_m, w_m  \rangle .$$
Set $T=\langle z_{t+1}, z_{t+3}, \cdots, z_{m-1}, w_{t+2}, w_{t+4}, \cdots, w_m\rangle$.
Noting that any two elements in $T$ are mutually exchanged and unequal involutions,
it deduces that
$T=\langle z_{t+1}\rangle \times \langle w_{t+2}\rangle \times  \cdots \times \langle z_{m-1}\rangle \times \langle w_{m}\rangle$ and $|T|=2^{m-t}$.
Since $z_{t+2s-1}$ and $w_{t+2s}$ are not contained in $H$ for $1\leqslant s\leqslant (m-t)/2$, we conclude that $T\cap H=1$.

Noticing that $H=\langle x_1,x_2,\cdots,x_m \rangle\Z(G)$ is abelian, $y_i\not\in H$ and $\Z(G)\leqslant \langle x_i, y_i\rangle$ for $1\leqslant i \leqslant t$, we conclude that $H\cap \langle x_i, y_i\rangle=\langle x_i  \rangle\Z(G)$, and further combining with the fact that $\langle x_i, y_i \rangle\cong \D_8$, one yields that there must exist an involution denoted by $v_i$ such that $v_i\in \langle x_i, y_i \rangle \setminus (H\cap \langle x_i, y_i  \rangle)$ and $v_ix_i\neq x_iv_i$.
Let $V=\langle v_i: 1\leqslant i\leqslant t\rangle$.
Then $v_iv_j=v_jv_i$ and so $V=\langle v_1 \rangle \times \langle v_2 \rangle \times \cdots \times \langle v_t \rangle$ and $|V|=2^t$.
On the other hand, the fact that $\langle x_1,y_1 \rangle \circ \cdots \circ  \langle x_t,y_t \rangle  \leqslant \C_G(T)$ leads to that $V \leqslant \C_G(T)$ and so $V T \subseteq \Omega_1(G)$.
However, $(\langle x_1,y_1 \rangle \circ \cdots \circ  \langle x_t,y_t \rangle) \cap T=1$ leads to that $V \cap T=1$ and so $VT=V\times T$ and $|V T|=2^m$.
We assert that $ H\cap VT=1$ and so $|\Omega_1(G)H|\geqslant 2^{m}$. Assume that $1\neq g\in H\cap VT$. Then $g=v_1^{a_1}v_2^{a_2}\cdots v_t^{a_t}z_{t+1}^{a_{t+1}}w_{t+2}^{a_{t+2}}\cdots z_{m-1}^{a_{m-1}}w_{m}^{a_m}$ for $a_i=0,1$ and $i=1,2,\cdots,m$.
Since $g\neq 1$, there exists $i$ such that $a_i=1$, and so without loss of generality, we can set $i=1$, i.e., $a_1=1$. Let $b=v_2^{a_2}\cdots v_t^{a_t}z_{t+1}^{a_{t+1}}w_{t+2}^{a_{t+2}}\cdots z_{m-1}^{a_{m-1}}w_{m}^{a_m}$. Then $g=v_1b$ and $bx_1=x_1b$.
However, $H$ is abelian, which leads to  $x_1g=gx_1$ and
$x_1v_1b=x_1g=gx_1=v_1bx_1=v_1x_1b$. Then $x_1v_1=v_1x_1$, a contradiction arisen, and hence the assertion holds.
It follows from $G/H\cong(G/\Z(G))/(H/\Z(G))$ is elementary abelian 2-group that $\Omega_1(G/H)=G/H$ and $|\Omega_1(G/H)|=|G|/|H|=2^{2m+1}/2^{m+1}=2^m$.
However, since $\Omega_1(G)H\subseteq\Omega_1(G/H)$ and $|\Omega_1(G)H|\geqslant 2^{m}$, one yields that $\Omega_1(G/H)=\Omega_1(G)H$. Employing the result of the Theorem \ref{theorem1.1}, we can deduce that $H$ is a perfect code of $G$. Hence the Claim 1 holds.

%
\textbf{Claim 2:} If $H$ is a perfect code of $G$, then $G\cong G_{m,1}$.

Assume that  $H$ is  a perfect code of $G$.  According to the case (2) of Lemma \ref{lemma2.4}, $G$ is either isomorphic to $G_{m,1}$ or $G_{m,2}$. Next, we will prove the claim 2 by reducing to absurdity. Suppose that $G\cong G_{m,2}$.
Applying the Lemma \ref{lemma2.5}, one yields that $m-t$ is odd.
Further, the case (3) of  Lemma \ref{lemma2.4} shows that $H\cong \ZZ_4\times \ZZ_2^{m-1}$ and $\exp H=4$.
Since $\Z(G)\cong \ZZ_2$ and $H=\langle x_1,x_2,\cdots, x_m\rangle \Z(G)$, there exists $x_i$ such that $o(x_i)=4$.
 Leaded from the case (4) of Lemma \ref{lemma2.4}, we deduce that $\Z(G)\leqslant \langle x_1,x_2,\cdots, x_m\rangle $ and $H=\langle x_1,x_2,\cdots, x_m\rangle$.
Let $z_j=y_{j}x_{j+1}$ and $w_{j+1}=y_{j+1}x_j$ for $t+1\leqslant j\leqslant m-2$.
By the complete similar argument to the  proof of the claim 1, we conclude that $\langle x_j,z_j \rangle\cong \langle x_j,w_{j+1} \rangle \cong \D_8$ for $t+1\leqslant j\leqslant m-2$ and $$G=\langle x_1, y_1  \rangle \circ \langle x_2, y_2  \rangle \circ \cdots  \circ  \langle x_t, y_t  \rangle \circ \langle x_{t+1}, z_{t+1}  \rangle \circ \langle x_{t+2}, w_{j+2}  \rangle \circ \cdots  \circ \langle x_{m-2}, z_{m-2}   \rangle  \circ \langle x_{m-1}, w_{m-1}  \rangle \circ \langle x_{m}, y_{m}  \rangle.$$
%
%
Set $I=\{i:o(x_i)=4, 1\leqslant i\leqslant t \}$,  
$J=\{j: (2,j-t)=1, o(x_j)=4,t+1\leqslant j\leqslant m-2 \}$ and  $K=\{k: (2,k-t)=2,o(x_k)=4, t+1\leqslant k\leqslant m-2 \}$. Let $g=\prod\limits_{i\in I}y_i   \prod\limits_{j\in J} z_{j}    \prod\limits_{k\in K} w_{k}  y_m$. According to the case (4) of Lemma \ref{lemma2.4}, $G$ has a unique square element $c$ is contained in $Z(G)\leqslant H$, it follows that $g^{2}\in H$.
%
%
%
Note that $H=\langle x_1,x_2,\cdots, x_m\rangle$. For any $h\in H$, we can set $h=x_1^{a_1}x_2^{a_2}\cdots x_m^{a_m}$, where $a_s$ are non-negative integers and $1\leqslant s\leqslant m$.
Since $x_iz_j=z_jx_i$ and $x_{i+1}w_{j+1}=w_{j+1}x_{i+1}$ for $i\neq j$,
we deduce that
$$hg=\prod\limits_{i\in I}x_i^{a_i}y_l    \prod\limits_{j\in J} x_{j}^{a_j}z_j    \prod\limits_{k\in K} x_k^{a_k}w_k     x_m^{a_m}y_m.$$
Simple calculation shows that $\prod\limits_{i\in I}x_iy_i$, $\prod\limits_{j\in J} x_{j}^{a_j}z_j$  and $\prod\limits_{k\in K} x_k^{a_k}w_k$ are involutions.
Noting that $\langle x_m,y_m\rangle\cong \Q_8$ and $o(x_m)=o(y_m)=4$, one yields that  $o(x_m^{a_m}y_m)=4$, and hence $hg$ is of order $4$.
By the arbitrariness of $h$, we know that $Hg$ has no involution.
Derived from $\Omega_1(G)H\subseteq\Omega_1(G/H)$ and $g^2\in H$, one yields that $Hg\in\Omega_1(G/H)\setminus\Omega_1(G)H$ and $\Omega_1(G)H$ is a proper subset of $\Omega_1(G/H)$. Employing the result of the Theorem \ref{theorem1.1}, we can deduce that $H$ is not a perfect code of $G$, a contradiction. Then the claim 2 holds and the Lemma \ref{maximal abelian subgroup} is right.           \qed

\noindent {\bf The proof of Theorem \ref{theorem1.2}.}

 \pf Suppose that $G$ is an extraspecial $2$-group and $H\leqslant G$.

($\Leftarrow$):  Assume that  either $H$ is non-abelian, or $H$ is abelian with $H\ntrianglelefteq G$, or $H$ is a maximal abelian subgroup of $ G\cong G_{m, 1}$.
Employing the result of the Lemma \ref{non-abelian} whenever $H$ is non-abelian, the Lemma \ref{nonormal} whenever $H$ is a normal non-maximal abelian and the Lemma \ref{maximal abelian subgroup} whenever $H$ is a maximal abelian, we can conclude that $H$ is a perfect code of $G$.

($\Rightarrow$): Suppose that $H$ is a perfect code of $G$.
Applying the Lemma \ref{normal non-maximal}, we deduce that $H$ is either  non-abelian, or  $H$ is abelian and not normal in $G$, or $H$ is a maximal abelian subgroup of $G$.
Further, if $H$ is a maximal abelian subgroup of $G$, then the Lemma \ref{maximal abelian subgroup} shows that $G\cong G_{m,1}$.

To sum up, we have proved Theorem \ref{theorem1.2}. \qed

\section{Proof of Theorem \ref{theorem1.3}}
In this section, we will prove the correctness of the Theorem \ref{theorem1.3}.
Let $G$ be a group with Sylow $2$-subgroup being extraspecial and $H\leqslant G$.
Then $H_2\leqslant\N_G(H_2)_2\leqslant P$ for $P$ a Sylow $2$-subgroup of $G$,
in particular, $P$ is an extraspecial $2$-group, i.e., $\Z(P)\cong\Z_2$, $P/\Z(P)$ is elementary abelian, in particular, $\exp P=4$.

\pf
($\Leftarrow$): If the order of $H$ is odd, then $H_2=1$  and the result of the Theorem \ref{theorem1.1} shows that $H$ is a perfect code of $G$. Thus, we assume that the order of $H$ is even as follows.
Assume that the case (1) or (2) of the Theorem \ref{theorem1.3} holds, that is $H_2$ is non-abelian or $H_2$ is abelian and $|\N_G(H_2)_2|<|P|$(i.e., $H_2 \ntrianglelefteq P$)  respectively.
Since $P$ is an extraspecial $2$-group,
according to Theorem \ref{theorem1.2},  $H_2$ is a perfect code of $P$ and so $H$ is a perfect code of $G$ by  Theorem \ref{theorem1.1}.
Suppose that the case (3) of the Theorem \ref{theorem1.3} holds, that is $P\cong G_{m,1}$ and $H_2$ is abelian with $|H_2|^2=2|P|=2|G_{m,1}|$, i.e., $|H_2|=2^{m+1}$.
The cases (3) of the Lemma \ref{lemma2.4} shows that $H_2$ is a maximal abelian subgroup of $P$.
Further, applying the Theorem \ref{theorem1.2},  $H_2$ is a perfect code of $P$ and so $H$ is a perfect code of $G$ following from the result of the Theorem \ref{theorem1.1}.

($\Rightarrow$): Suppose that $H$ is a  perfect code of $G$.
According to Theorem \ref{theorem1.1}, we conclude that $H_2$ is a perfect code of $P$. Now assume that $|H|$ is even, that is $H_2\neq 1$.
%
Employing the result of Theorem \ref{theorem1.2}, we deduce that $H_2$ is either non-abelian, or abelian and $H\ntrianglelefteq P$, or maximal abelian subgroup of $P\cong G_{m,1}$.
Assume that $H_2$ is a  maximal abelian subgroup of $P$.  Applying the case (3) of Lemma \ref{lemma2.4}, one yields that $|H_2|^{2}=2|P|$.

To sum up, we have proved Theorem \ref{theorem1.3}. \qed

 \noindent{\bf Declaration of competing interest}

The authors declare that they have no known competing financial interests or personal relationships that could have appeared to influence the work reported in this paper.

 \noindent{\bf Acknowledgments}

%
The authors would like to thank Prof. Binzhou Xia and Prof. Junyang Zhang  for their helpful suggestions, which greatly improved the quality of this work.


\begin{thebibliography}{99}



\bibitem{c}
Caliskan,C., Miklavi$\mathrm{\check{c}}$,$\mathrm{\check{S}}$., $\mathrm{\ddot{O}}$zkan,S.,  $\mathrm{\check{S}}$parl,P.: Efficient domination in Cayley graphs of generalized dihedral groups. Discuss. Math. Graph Theory. 42, 823-841(2022).

\bibitem{ch}
Chen,J.,  Wang,Y.,   Xia,B.: Characterization of subgroup perfect codes in Cayley graphs. Discrete Math.  343, 111813(2020).

\bibitem{dej}
Dejter,I.J., Serra,O.: Efficient dominating sets in Cayley graphs. Discrete Appl. Math. 129,319-328(2003).

\bibitem{den}
 Deng,Y.-P., Sun,Y.-Q., Liu,Q., Wang,H.-C. :Efficient dominating sets in circulant graphs. Discrete Math.  340, 1503-1507(2017).

\bibitem{f}
R. Feng,  H. Huang, S. Zhou,  Perfect codes in circulant graphs, Discrete Math. 340(7)(2017) 1522-1527.

\bibitem{h}
Huang,H., Xia,B.,  Zhou,S.: Perfect codes in Cayley graphs. SIAM J. Discrete Math. 32, 548-559(2018).

\bibitem{hu}
Huppert, B.: {Endliche Gruppen I}. Springer-Verlag, Berlin, 1967.

\bibitem{j}
Kratochv$\acute{i}$l,J.: Perfect codes over graphs. J. Combin. Theory, Ser B. 40, 224-228(1986).

\bibitem{kh}
 Khaefi,Y.,  Akhlaghi,Z.,  Khosravi,B. : On the subgroup perfect codes in Cayley graphs, Des. Codes  Cryptogr. 91, 55-61(2023).

\bibitem{ku}
Kumar,K.R.,  MacGillivray,G.: Efficient domination in circulant graphs. Discrete Math.  313,767-771(2013).

\bibitem{l}
 Lee,J.: Independent perfect domination sets in Cayley graphs. J. Graph Theory.  37, 213-219(2001).


\bibitem{m}
Ma,X., Walls,G.L., Wang,K., Zhou,S.: Subgroup perfect codes in Cayley graphs. SIAM J. Discrete Math. 34, 1909-1921(2020).


\bibitem{o}
Obradovi$\mathrm{\acute{c}}$,N.   Peters,J.,  Ru$\mathrm{\check{z}}$i$\mathrm{\acute{c}}$,G.:   Efficient domination in circulant graphs with two chord lengths. Inf. Process. Lett. 102, 253-258(2007).


\bibitem{z}
Zhang,J.: Characterizing subgroup perfect codes by 2-subgroups.  Des. Codes Cryptogr. 91,2811-2819(2023).

\bibitem{z2}
Zhang,J., Zhou,S.: On subgroup perfect codes in Cayley graphs. Eur. J. Comb.  91,103228(2021).

\bibitem{z3}
Zhang, J., Zhou, S.:  Corrigendum to ``On subgroup perfect codes in Cayley graphs [Eur. J. Comb. 91,103228].".  Eur. J. Comb. 101,103461(2022).

\end{thebibliography}
\end{document}